\documentclass[11pt]{article}

\date{{November 9, 2026}}

\PassOptionsToPackage{dvipsnames,svgnames,x11names}{xcolor}
\usepackage{xcolor}
\definecolor{labelkey}{rgb}{0,0.08,0.45}
\definecolor{refkey}{rgb}{0,0.6,0.0}
\definecolor{Brown}{rgb}{0.45,0.0,0.05}
\definecolor{lime}{rgb}{0.00,0.8,0.0}
\definecolor{lblue}{rgb}{0.5,0.5,0.99}
\definecolor{OliveGreen}{rgb}{0,0.6,0}
\definecolor{tyrianpurple}{rgb}{0.4, 0.01, 0.24}
\definecolor{myseagreen}{HTML}{3FBC9D}
\definecolor{myblue}{rgb}{0.9,0.9,0.98}
\colorlet{hlcyan}{cyan!30}

\newif\ifshowdraftlabels
\showdraftlabelsfalse
\ifshowdraftlabels\usepackage{showkeys}\fi

\usepackage{lmodern} 
\usepackage{mathpazo} 
\usepackage{microtype}
\usepackage{subcaption} 
\usepackage[T1]{fontenc}
\usepackage{mathtools}
\usepackage{amssymb}
\IfFileExists{stmaryrd.sty}{\usepackage{stmaryrd}}{}
\usepackage{empheq}

\usepackage{hyperref}
\hypersetup{
  colorlinks=true,
  linkcolor=refkey,
  urlcolor=lblue,
  citecolor=magenta
}

\usepackage{amsthm}
\makeatletter
\def\th@plain{%
  \thm@notefont{}%
  \itshape 
}
\def\th@definition{%
  \thm@notefont{}%
  \normalfont 
}

\usepackage[capitalize,nameinlink]{cleveref}
\makeatother
\newtheorem{theorem}{Theorem}[section]
\newtheorem{lemma}[theorem]{Lemma}
\newtheorem{corollary}[theorem]{Corollary}
\newtheorem{proposition}[theorem]{Proposition}

\newtheorem{example}[theorem]{Example}

\newtheorem{remark}[theorem]{Remark}
\crefname{theorem}{Theorem}{Theorems}
\Crefname{theorem}{Theorem}{Theorems}
\crefname{fact}{Fact}{facts}
\Crefname{fact}{Fact}{facts}
\crefname{equation}{}{equations}
\crefname{chapter}{Appendix}{chapters}
\crefname{item}{}{items}
\crefname{enumi}{}{}

\IfFileExists{nicematrix.sty}{\usepackage{nicematrix}}{}
\usepackage{booktabs}
\IfFileExists{siunitx.sty}{\usepackage{siunitx}}{}
\usepackage{tikz}
\usepackage[color=myseagreen]{todonotes}
\IfFileExists{soul.sty}{\usepackage{soul}}{}
\usepackage{nameref}
\usepackage{graphicx}
\usepackage{float}
\usepackage[shortlabels,inline]{enumitem}
\setlist[enumerate]{nosep}
\usepackage{relsize}

\usepackage[margin=1in,footskip=0.4in]{geometry}
\makeatletter

\let\orig@label\label

\renewcommand{\label}[1]{%
  \begingroup
  \def\@currentlabelname{}%
  \ifx\current@theorem\relax\else
    \def\@currentlabelname{\current@theorem}%
  \fi
  \ifx\cref@currentlabel\undefined\else
    \let\@currentlabelname\cref@currentlabel
  \fi
  \orig@label{#1}%
  \endgroup
}

\AtBeginEnvironment{theorem}{\def\current@theorem{theorem}}
\AtBeginEnvironment{proposition}{\def\current@theorem{proposition}}
\AtBeginEnvironment{fact}{\def\current@theorem{fact}}
\makeatother

\newcommand{\seppthree}{\setlength{\itemsep}{-3pt}}

\newcommand{\Pj}[1]{\ensuremath{\operatorname{P}_{#1}}}

\newcommand{\bId}{\ensuremath{\operatorname{\mathbf{Id}}}}

\newcommand{\bx}{\ensuremath{\mathbf{x}}}

\newcommand{\bzero}{\ensuremath{\boldsymbol{0}}}

\newcommand{\weakly}{\ensuremath{\:{\rightharpoonup}\:}}

\newcommand{\thalb}{\ensuremath{\tfrac{1}{2}}}
\newcommand{\menge}[2]{\big\{{#1}~\big|~{#2}\big\}}

\newcommand{\scal}[2]{\left\langle{#1},{#2}\right\rangle}

\newcommand{\tscal}[2]{\langle{#1},{#2}\rangle}

\newcommand{\RR}{\ensuremath{\mathbb{R}}}
\newcommand{\HH}{\ensuremath{\mathcal{H}}}

\newcommand{\RP}{\ensuremath{\mathbb{R}_+}}
\newcommand{\RPP}{\ensuremath{\mathbb{R}_{++}}}

\DeclareMathOperator*{\argmin}{argmin}

\newcommand{\prox}{\ensuremath{\operatorname{P}}}

\newcommand{\cspan}{\ensuremath{\overline{\operatorname{span}}\,}}

\newcommand{\Fix}{\ensuremath{\operatorname{Fix}}}
\newcommand{\Id}{\ensuremath{\operatorname{Id}}}

\newcommand{\bX}{\ensuremath{\mathbf{X}}}
\newcommand{\bH}{\ensuremath{\mathbf{H}}}

\newcommand{\bK}{\ensuremath{\mathbf{K}}}

\newcommand{\by}{\ensuremath{\mathbf{y}}}
\newcommand{\bz}{\ensuremath{\mathbf{z}}}

\providecommand{\norm}[1]{\lVert#1\rVert}

\providecommand{\bw}{\mathbf{w}}

\newcommand{\mybluebox}[1]{\colorbox{myblue}{\hspace{1em}#1\hspace{1em}}}

\author{
Heinz H.\ Bauschke\thanks{\raggedright
Mathematics, University
of British Columbia,
Kelowna, B.C.\ V1V~1V7, Canada. E-mail:
\texttt{heinz.bauschke@ubc.ca}.}~~~and~
Tran Thanh Tung\thanks{\raggedright
Mathematics, University
of British Columbia,
Vancouver, B.C.\ V6T~1Z4, Canada. 
 E-mail: \texttt{tung.tran@ubc.ca}.}
}

\title{\textsf{{
A family of alternating projection examples\\ whose 
Ces\`aro means fail to converge strongly}}}

\begin{document}
\allowdisplaybreaks
\maketitle
\enlargethispage{\baselineskip}

\begin{abstract}
In 1965, Bregman proved that the sequence of iterates of the method of alternating projections 
applied to two closed convex sets with nonempty intersection in a Hilbert space 
converges weakly to a point in the intersection. The question whether 
this convergence may actually fail to be strong was settled in 
a seminal paper by Hundal in 2004. 
He provided a counterexample featuring a halfspace and 
a closed convex cone generated by a certain curve. 
While largely elementary, the details of Hundal's construction and analysis are quite involved. 

We revisit this seminal example and present a family of such examples that 
allows for simpler proofs. 
Our generating curve combines unit vectors with Gaussian inner products
and a height determined by a scalar function.
The orbits converge weakly --- but not strongly --- to the origin, while their norms have a positive
limit for which we give explicit lower bounds. We also prove 
that the associated Ces\`aro means behave likewise; in fact, this 
is generalized to alternating proximal mappings whose composition 
has a fixed point. 
Finally, these orbits allow interpretations in the context of 
projected-gradient, proximal-point, and gradient-descent algorithms. 
\end{abstract}



\small
\section{Introduction}
\label{s:intro}

Throughout this paper, 
we assume that 
\begin{empheq}[box=\mybluebox]{equation}
  \label{e:X}
  \text{$X$ is an \emph{infinite-dimensional} real Hilbert space,}
\end{empheq}
with inner product $\scal{\cdot}{\cdot}$ and induced norm
$\norm{\cdot}$.

The method of alternating projections is one of the most 
fundamental algorithms for
finding a point in the nonempty intersection of two closed convex sets. 
In 1965, Bregman
\cite{Bregman1965} proved that the
iterates converge \emph{weakly} to a point in the intersection. 
The question whether strong convergence always holds 
--- as it does for two closed linear subspaces ---
remained open for many years. In his seminal 2004 paper, 
Hundal \cite{Hundal} showed that strong convergence may fail for
general closed convex sets. His example uses a closed convex cone and a
halfspace; along the constructed orbit, the halfspace can be replaced by
its boundary hyperplane. The intersection consists only of the origin.
The iterates converge weakly to the origin, while their norms have a
strictly positive limit.

Modifications and simplifications of Hundal's construction were given by
Matou\v{s}kov\'a and Reich \cite{MR2003} and by Kopeck\'a
\cite{Kopecka2008}. 
Their constructions also cover disjoint sets, and
Kopeck\'a obtains a cone with nonempty interior. 
Related examples appear also 
in fixed point theory and convex optimization: Genel and Lindenstrauss
\cite{GL1975} constructed weakly but not strongly convergent iterates of
a firmly nonexpansive mapping, while G\"uler \cite{Guler1991} obtained this
behaviour even for the proximal point algorithm. 
In \cite{BMR}, it is shown  that Hundal's example also yields
counterexamples for averaged projections and proximal mappings.

It is natural to ask whether averaging improves convergence. 
Baillon's
nonlinear mean ergodic theorem \cite{Baillon1975} guarantees weak
convergence of the Ces\`aro means of a nonexpansive mapping with fixed
points. 
Strong convergence need not follow --- this was already discussed
by Baillon \cite[p.~VIII.29]{BaillonSeminar1978}. 
Krengel and Lin
\cite[Theorem~5.1 and Example~5.2]{KL1987} give an explicit example in
$L^2$ in which even the norms of the means fail to converge. 
Gaussian
constructions for firmly nonexpansive mappings were recently 
given in \cite{BT26}, and already utilized 
by Quan and Patrinos \cite{QP2026} to establish failure 
of strong convergence for the Douglas-Rachford splitting algorithm in 
the context of convex feasibility. 

For alternating projections onto two closed convex sets with nonempty
intersection, the averaging question has a very useful additional feature:
each of the two projection subsequences is asymptotic in norm to its own Ces\`aro
means. We give a proof of a more general result in 
\cref{t:Phi}. 
Consequently, Hundal's original example
already implies failure of strong convergence of the Ces\`aro means!
\emph{
 Our goal in this paper is twofold. 
 First, we show that strong cluster points of alternating projections 
 are precisely those of their Ces\`aro means. Second, we
 construct a new family of cones for which
 alternating projections with a
 fixed hyperplane carry each generating ray to the next, and we 
 show that every nonzero orbit so obtained converges weakly but not
 strongly to zero. It appears that this construction is significantly simpler 
 than the ones previously available in the literature.
}

For the rest of this paper, we assume that we are given a Gaussian curve of unit
vectors that satisfies
\begin{empheq}[box=\mybluebox]{equation}
  \label{e:Kurve}
  u\colon\RP\to X
  \quad\text{and}\quad
  (\forall s,t\in\RP)\quad \scal{u(s)}{u(t)}=e^{-(s-t)^2}.
\end{empheq}
Such curves exist\footnote{For instance, if $(e_k)_{k\geq0}$ is an orthonormal
sequence in $X$, then one may take
\begin{equation*}
  u(t):=e^{-t^2}\sum_{k=0}^{\infty}
                \sqrt{\frac{2^k}{k!}}\,t^k e_k.
\end{equation*}} in every infinite-dimensional real Hilbert space 
\cite{BT26}; 
moreover, the curve is differentiable and
\begin{equation}
  \label{e:Kweak}
  u(t)\weakly0\quad\text{as }t\to\infty.
\end{equation}
We will also assume that
\begin{subequations}
  \label{e:f}
  \begin{empheq}[box=\mybluebox]{align}
    \label{e:f0}
    f
    \;\text{is twice continuously differentiable on
    an open interval containing $\RP$,}
  \end{empheq}
and we impose the following conditions\footnote{Note that \cref{e:f0}
and \cref{e:f3} make the integral in \cref{e:f2} well defined.} on $f$:
  \begin{empheq}[box=\mybluebox]{align}
    f'(0)e^{-2f(0)}&\leq\thalb,\label{e:f1}\\[0.2em]
    \int_0^\infty\frac{dt}{f'(t)}&<\infty,\label{e:f2}\\[0.2em]
    (\forall t\in\RP)\quad 2&\leq f'(t)-2t,\label{e:f3}\\[0.2em]
    (\forall t\in\RP)\quad 0&\leq f''(t),\label{e:f4}\\[0.2em]
    (\forall t\in\RP)\quad f''(t)&\leq
        2+\thalb\big(f'(t)-2t\big)^2.\label{e:f5}
  \end{empheq}
\end{subequations}
These assumptions admit simple examples, including $f(t)=(1+t)^3$
and $f(t)=2e^t$  (see \cref{ex:powers} and \cref{ex:exponential} below 
as well as \cref{s:examples} for a convenient sufficient condition).

The two sets we project upon reside in the real Hilbert space
\begin{empheq}[box=\mybluebox]{align}
  \label{e:bX}
  \bX:=\RR\times X,
\end{empheq}
with the usual inner product and norm. 
We now define the hyperplane
\begin{empheq}[box=\mybluebox]{align}
  \label{e:bH}
  \bH:=\{0\}\times X,
\end{empheq}
the curve
\begin{empheq}[box=\mybluebox]{align}
  \label{e:bw}
  \bw\colon\RP\to\bX\colon t\mapsto\big(e^{-f(t)},u(t)\big),
\end{empheq}
and the cone
\begin{empheq}[box=\mybluebox]{align}
  \label{e:bK}
  \bK:=\text{the closed convex cone generated by }
              \{\bw(t)\}_{t\in\RP}.
\end{empheq}
We will show in \cref{p:HK0} that $\bH\cap\bK=\{\bzero\}$.
Given $s_0\geq0$ and $\rho_0>0$, set $\by_0:=\rho_0\bw(s_0)$ and
generate the alternating projections by
\begin{equation*}
  (\forall n\geq0)\quad
  \bx_{n+1}:=\Pj{\bH}\by_n
  \quad\text{and}\quad
  \by_{n+1}:=\Pj{\bK}\bx_{n+1}.
\end{equation*}
Both sequences converge weakly to $\bzero$, while their norms converge
to the same number $\rho_\infty>0$ (see \cref{t:map} and \cref{p:rhon}). 
A general averaging result then implies that the
Ces\`aro means of both subsequences converge weakly to $\bzero$ and have
the same positive limiting norm $\rho_\infty$ (see \cref{t:cesaro}).
A more detailed comparison with
Hundal's construction is given in \cref{r:Hundal}.
Finally, the standard relations between projections and convex
optimization give projected-gradient, proximal-point, and
gradient-descent interpretations of these examples. For the latter
two interpretations, we use the average of the projections onto
$\bK$ and its reflection across $\bH$, as in \cite{BMR}.
These interpretations are discussed in \cref{r:optimization}.

The remainder of this paper is organized as follows.
In \cref{s:aux}, we provide examples of $f$ and establish the inequalities
used in the construction. The interplay between alternating 
proximal mappings and their Ces\`aro means is studied in \cref{s:Phi}. 
\cref{s:sets} describes the geometry of $\bH$
and $\bK$ and gives the required projection formulas. In \cref{s:seq},
we study the scalar sequences determining the alternating projections
and prove the positive lower bounds for their limiting radius.
\cref{s:main} contains the weak-but-not-strong convergence result, 
and resulting interpretations for convex optimization algorithms. 

\section{Auxiliary results}
\label{s:aux}

\subsection{\texorpdfstring{Incarnations of $f$}{Incarnations of f}}
\label{s:examples}

We first give explicit examples of functions satisfying \cref{e:f}.
We then provide a simple sufficient condition that is preserved
under convex combinations.

\begin{example}
\label{ex:powers}
  If
  \begin{equation}
    f\colon \left]-1,\infty\right[\to\RR\colon 
    t\mapsto (1+t)^p, \quad\text{with}\quad 2<p\leq 1+\sqrt{5},
  \end{equation}
  then $f$ satisfies \cref{e:f}.
\end{example}
\begin{proof}
  We have $f'(t) = p(1+t)^{p-1}$ and $f''(t) = p(p-1)(1+t)^{p-2}$,
  the latter of which is continuous and nonnegative.
  Hence \cref{e:f0} and \cref{e:f4} hold.

  Next, we have ${\color{black}f'(t)-2t} = p(1+t)^{p-1}-2t$,
  and ${\color{black}e^{-2f(t)}} = {\color{black}e^{-2(1+t)^p}}$. Then $f'(0){\color{black}e^{-2f(0)}} = p{\color{black}e^{-2}}\leq (1+\sqrt{5})/{\color{black}e^{2}}
    \approx 0.438 < \thalb$, so \cref{e:f1} holds.

  To verify \cref{e:f2}, we compute
  \begin{equation*}
    \int_0^\infty \frac{dt}{f'(t)} = \int_0^\infty \frac{dt}{p(1+t)^{p-1}}
    = \frac{1}{p}\int_1^\infty u^{1-p} du = 
    \frac{1}{p}\frac{1}{2-p} u^{2-p}\Big|_1^\infty = 
    \frac{1}{p(p-2)} < \infty
  \end{equation*}
  and note that $p>2$ is needed for the integral to be finite.

  Next, ${\color{black}f'(t)-2t} = p(1+t)^{p-1}-2t \geq 2$ is equivalent to
  $p(1+t)^{p-1} \geq 2t+2 = 2(1+t)$,
  which is true because $p\geq 2$. This yields \cref{e:f3}.

  Finally, let us verify \cref{e:f5}.
  Let $t\in\RP$. Because $p\geq 2$, we have
  $s := (1+t)^{p-2}\geq 1$, $f''(t)=p(p-1)s$, and
  $f'(t)-2t= p(1+t)^{p-1}-2t = ps(1+t)-2t = ps+ t(ps-2) \geq ps$.
  It thus suffices to show that
    $p(p-1)s \leq 2+\thalb (ps)^2$, i.e., 
    $2+\thalb(ps)^2 - p(p-1)s\geq 0$. 
  Now rewrite $2+\thalb (ps)^2 - p(p-1)s$
  as a sum of three terms:
  \begin{equation*}
    \big( 2+\thalb p^2 - p(p-1)\big)  + p(s-1)
    + \thalb p^2(s-1)^2.
  \end{equation*}
  The last two summands are nonnegative, so it suffices to
  show that
  \begin{equation*}
    2+\thalb p^2 - p(p-1) \geq 0.
  \end{equation*}
  In turn, this quadratic inequality is equivalent
  to $p^2 - 2p -4 \leq 0$, and to
  $p \in [1-\sqrt{5}, 1+\sqrt{5}]$.
  Because $2<p\leq 1+\sqrt{5}$, this holds, and so \cref{e:f5} is verified.
\end{proof}

One may construct further examples of $f$ where \cref{e:f} holds --- we 
present one sufficient condition that is fairly easy to check and also closed 
under convex combinations and multiplication by scalars greater 
than or equal to one: 

\begin{proposition}
\label{p:cutesuff}
Let $f$ be twice continuously differentiable on 
an open interval containing $\RP$ and suppose that 
\begin{equation}
\label{e:cutesuff}
f'(0)\leq 2f(0) \quad\text{and}\quad (\forall t\in\RP)\quad 
t+2 \leq f''(t)\leq f'(t). 
\end{equation}
Then $f$ satisfies \cref{e:f}.
\end{proposition}
\begin{proof}
Note that $2\leq f''(0)\leq f'(0)$ and so 
$1 \leq \thalb f'(0)\leq f(0)$. 
Integrating $t+2\leq f''(t)$, we thus obtain 
\begin{equation}
\label{e:261103f}
\thalb(t+2)^2 = 2+2t+\thalb t^2 \leq f'(0)+2t +\thalb t^2\leq f'(t). 
\end{equation}
Because $f'(0)\leq 2f(0)$, $1\leq f(0)$, and 
$\left[\thalb,\infty\right[\to\RR \colon 
t\mapsto t e^{-2t}$ is strictly decreasing, we have
$f'(0)e^{-2f(0)}\leq 2f(0)e^{-2f(0)}\leq 
2e^{-2}<\thalb$ and so \cref{e:f1} holds.
Next, we deduce from \cref{e:261103f} that
\begin{equation*}
\int_0^\infty \frac{dt}{f'(t)} \leq
\int_0^\infty \frac{2dt}{(t+2)^2}
= \frac{-2}{t+2}\Big|_0^\infty = 1 < \infty, 
\end{equation*}
so \cref{e:f2} holds. 
\cref{e:261103f} also implies that $f'(t)-2t\geq 2+\thalb t^2$ which yields \cref{e:f3}.
Clearly, \cref{e:f4} follows from \cref{e:cutesuff}. 

Now set $h(t) := f'(t)-2t$. 
Then \cref{e:261103f} yields 
$h(t)\geq 2 +\thalb t^2 \geq 2t$. 
Therefore, using \cref{e:cutesuff}, we estimate 
\begin{align*}
f''(t) 
&\leq f'(t) = h(t)+2t
\leq 2h(t) 
\leq 2+ \thalb\big(h(t)\big)^2
= 2+\thalb\big(f'(t)-2t\big)^2,
\end{align*}
and this establishes \cref{e:f5}.
\end{proof}

\begin{example}
\label{ex:exponential}
  If $\alpha\in\left[2,\infty\right[$ and 
  \begin{equation}
    f = \alpha\exp, 
  \end{equation}
  then $f$ satisfies \cref{e:f}.
\end{example}
\begin{proof}
Indeed, the assumptions of \cref{p:cutesuff} hold. 
\end{proof}

\subsection{\texorpdfstring{Properties of $f$}{Properties of f}}

We now return to a general function $f$ satisfying \cref{e:f}
and establish some useful properties.

\begin{lemma}
  \label{ff}
  Recall that \cref{e:f} holds. Then the following hold:
  \begin{enumerate}[label=(\roman*), ref=(\roman*), font=\normalfont]
    \item \label{ff:i} For every $t\in\RP$,
      \begin{equation*}
        f(t)>t^2
        \quad\text{and}\quad
        f'(t)\geq 2(1+t)>0.
      \end{equation*}
      Consequently, $f(t)\to\infty$ and $e^{-2f(t)}\to0$ as $t\to\infty$.

    \item \label{ff:ii} For every $t\in\RP$,
      \begin{equation*}
        0\leq f''(t)\leq\big(f'(t)\big)^2.
      \end{equation*}
      Moreover, $1/f'$ is decreasing on $\RP$.

    \item \label{ff:iii} For every $t\in\RP$,
      \begin{equation*}
        0<e^{-2f(t)}\leq e^{-2f(0)}\leq\tfrac14
        \quad\text{and}\quad
        (e^{-2f})'(t)=-2f'(t)e^{-2f(t)}<0.
      \end{equation*}
      Consequently, $t\mapsto e^{-2f(t)}$ is strictly decreasing.
  \end{enumerate}
\end{lemma}
\begin{proof}
  \ref{ff:i}: By \cref{e:f3}, we have
  $f'(t)\geq2(1+t)>0$ for every $t\in\RP$.
  In particular, $f'(0)\geq2$. Combining this with \cref{e:f1}, we obtain
  \begin{equation}
    \label{e:261103a}
    e^{-2f(0)}\leq\frac{1}{2f'(0)}\leq\frac14.
  \end{equation}
  Hence $f(0)\geq\ln2>0$, and we obtain
  \begin{align*}
    f(t)
    &=f(0)+\int_0^t f'(\tau)\,d\tau
    \geq\ln2+\int_0^t2(1+\tau)\,d\tau
     =\ln2+2t+t^2>t^2.
  \end{align*}
  The same estimate implies that $f(t)\to\infty$, 
  and hence that $e^{-2f(t)}\to0$, as $t\to\infty$.

  \ref{ff:ii}: Let $t\in\RP$.
  By \cref{e:f3},
  $2\leq f'(t)-2t\leq f'(t)$. 
  This implies that 
  $4\leq (f'(t)-2t)^2$ and 
  $(f'(t)-2t)^2\leq (f'(t))^2$. 
  The former inequality is equivalent to 
  $2+\frac12(f'(t)-2t)^2\leq (f'(t)-2t)^2$. 
  Altogether, and with the help of 
  \cref{e:f4} and \cref{e:f5}, we estimate 
  \begin{align*}
    0\leq f''(t)
    &\leq2+\tfrac12\big(f'(t)-2t\big)^2
    \leq\big(f'(t)-2t\big)^2
     \leq\big(f'(t)\big)^2.
  \end{align*}
  Furthermore, using $f'>0$ from \ref{ff:i}, we obtain
  \begin{equation*}
    (1/f')'(t)=-\frac{f''(t)}{(f'(t))^2}\leq0.
  \end{equation*}
  Hence $1/f'$ is decreasing.

  \ref{ff:iii}: Differentiating and recalling \ref{ff:i}, we obtain
    $(e^{-2f})'(t)=-2f'(t)e^{-2f(t)}<0$.
  Thus $t\mapsto e^{-2f(t)}$ is strictly decreasing, and 
  \cref{e:261103a} implies 
  \begin{equation*}
    0<e^{-2f(t)}\leq e^{-2f(0)}\leq\tfrac14.
  \end{equation*}
  The proof is complete.
\end{proof}

\subsection{\texorpdfstring{Properties of $Q$}{Properties of Q}}

We now define an auxiliary function $Q$ and establish the
properties needed for our construction. 
Recall \cref{e:f}, and set 
\begin{empheq}[box=\mybluebox]{equation}
  \label{e:Q}
  Q\colon\RP\to\RR\colon t\mapsto
  \frac{f'(t){\color{black}e^{-2f(t)}}}{2(1+{\color{black}e^{-2f(t)}})} 
  = \frac{f'(t)}{2\big({\color{black}e^{2f(t)}}+1\big)}.
\end{empheq}

\begin{proposition}[properties of $Q$]
  \label{p:Q}
  The following hold:
  \begin{enumerate}[label=(\roman*), font=\normalfont]
    \item \label{p:Q1} $Q$ is continuous.
    \item \label{p:Q2} 
    For every $t\in\RP$, $Q(t)>0$.
    \item \label{p:Q3} $Q$ is strictly decreasing.
    \item \label{p:Q4} $Q(t)\to 0$ as $t\to\infty$.
  \end{enumerate}
  Now fix $s\in\RP$.
  Then there exists a unique $\color{black}s_+\in\RPP$ such that
  \begin{equation}
    \label{s+-s}
    s_+ - s = Q(s_+).
  \end{equation}
  Moreover, $s$ and $s_+$ satisfy the following inequalities:
\begin{equation}
    \label{svss+}
    s<s_+<s+Q(s)\leq s+Q(0)<s+\tfrac14,
  \end{equation}
  \begin{equation}
    \label{Es+}
    {\color{black}e^{-2f(s_+)}} < 
    \frac{5}{2} \int_{s}^{s_+} \frac{dt}{f'(t)},
  \end{equation}
  and
  \begin{equation}
    \label{Es-Es+}
    {\color{black}e^{-2f(s)}}-{\color{black}e^{-2f(s_+)}}> 4(s-s_+)^2.
  \end{equation}
\end{proposition}
\begin{proof}
  \cref{p:Q1}: Clear from \cref{e:f0} and \cref{e:Q}.

  \cref{p:Q2}: By \cref{ff}\cref{ff:i}\&\cref{ff:iii}, both $f'(t)$ and $e^{-2f(t)}$
  are positive, so $Q(t)>0$.

  \cref{p:Q3}: We work with $\ln\circ Q$.
  For every $t\in\RP$, we have
  \begin{equation*}
    (\ln\circ Q)(t)
    =\ln(f'(t))-2f(t)-\ln2-\ln\big(1+e^{-2f(t)}\big).
  \end{equation*}
  Using \cref{ff}, we obtain
  \begin{subequations}
  \label{e:261103b}
  \begin{align}
    (\ln\circ Q)'(t)
    &=\frac{f''(t)}{f'(t)}-2f'(t)
      +\frac{2f'(t)e^{-2f(t)}}{1+e^{-2f(t)}}
    =\frac{f''(t)}{f'(t)}-\frac{2f'(t)}{1+e^{-2f(t)}}\\
    &=f'(t)\left(\frac{f''(t)}{(f'(t))^2}
                  -\frac{2}{1+e^{-2f(t)}}\right)
    \leq f'(t)\left(1-\frac{2}{1+1/4}\right)
     =-\frac35 f'(t)\\
     &\leq-\frac65.
  \end{align}
\end{subequations}
  Hence $\ln\circ Q$ is strictly decreasing and thus so is $Q$. 

  \cref{p:Q4}:
  Integrating \cref{e:261103b} over $[0,t]$, we obtain
  \begin{equation*}
    \ln Q(t)-\ln Q(0)\leq-\frac65t,
  \end{equation*}
  and therefore
  \begin{equation*}
    0<Q(t)\leq Q(0)e^{-6t/5}\to0
    \quad\text{as }t\to\infty.
  \end{equation*}

  Now fix $s\in\RP$.
  By \cref{p:Q1}--\cref{p:Q3},
  the function $h\colon\RP\to\RR\colon t\mapsto t-s-Q(t)$
  is continuous and strictly increasing, with $h(s)=-Q(s)<0$.
  Also, since $Q(t)\leq Q(0)$,
  \begin{equation*}
    h(t)\geq t-s-Q(0)\to\infty
    \quad\text{as }t\to\infty.
  \end{equation*}
  The Intermediate Value Theorem therefore gives a unique
  $s_+>s$ such that $h(s_+)=0$, which is \cref{s+-s}.
  Since $Q$ is strictly decreasing,
  \begin{equation*}
    s_+=s+Q(s_+)<s+Q(s)\leq s+Q(0).
  \end{equation*}
  Finally, \cref{e:f1}  yields 
  \begin{equation*}
    Q(0)=\frac{f'(0)e^{-2f(0)}}{2(1+e^{-2f(0)})}
    \leq\frac{1}{4(1+e^{-2f(0)})}<\frac14.
  \end{equation*}
  We've thus verified  \cref{svss+}.

  Next, since $ 0 \leq s < s_+$, 
  \cref{ff}\ref{ff:iii} gives
    $e^{-2f(s_+)}<e^{-2f(0)}\leq\tfrac14$.
  This implies $4e^{-2f(s_+)}<1$, which is equivalent to 
  \begin{equation}
  \label{e:261103c}
    2\big(1+e^{-2f(s_+)}\big)< \frac{5}{2}. 
  \end{equation}
  By \cref{ff}\cref{ff:ii}, the function $1/f'$ is decreasing.
  Consequently, for every $t\in[s,s_+]$, we have
  $1/f'(t)\geq1/f'(s_+)$, and hence
  \begin{equation*}
    \frac{s_+-s}{f'(s_+)}\leq\int_s^{s_+}\frac{dt}{f'(t)}.
  \end{equation*}
  Recalling \cref{e:Q} and \cref{s+-s}, the positivity of $f'$ in
  \cref{ff}\ref{ff:i}, and \cref{e:261103c}, we therefore obtain
  \begin{align*}
    e^{-2f(s_+)}
    &=\frac{2(1+e^{-2f(s_+)})(s_+-s)}{f'(s_+)}
    <\frac52\cdot \frac{s_+-s}{f'(s_+)}
     \leq\frac52\int_s^{s_+}\frac{dt}{f'(t)}.
  \end{align*}
  Thus \cref{Es+} holds.

  Next, the derivative identity in \cref{ff}\cref{ff:iii}
  and the definition of $Q$ (see \cref{e:Q}) give
  \begin{equation}
    \label{e:261103d}
    -(e^{-2f})'(t)=2f'(t)e^{-2f(t)}
      =4(1+e^{-2f(t)})Q(t).
  \end{equation}
  For $t\in[s,s_+]$, we have $e^{-2f(t)}>0$ and
  $Q(t)\geq Q(s_+)>0$, by 
  \cref{p:Q2}, \cref{p:Q3}, and \cref{ff}\ref{ff:iii}. 
  Integrating \cref{e:261103d} over $[s,s_+]$, 
  and using \cref{s+-s}, we thus obtain
  \begin{align*}
    e^{-2f(s)}-e^{-2f(s_+)}
    &=4\int_s^{s_+}(1+e^{-2f(t)})Q(t)\,dt\\
    &>4\int_s^{s_+}Q(s_+)\,dt
    =4(s_+-s)Q(s_+)=4(s_+-s)^2.
  \end{align*}
  \nopagebreak[4]
  This verifies \cref{Es-Es+}, and we're done.
\end{proof}

For a given $s\in\RP$, \cref{p:Q} determines $s_+$ uniquely (although
we do not have a closed-form expression for it). 
We therefore now well-define 
the operator\footnote{It can be shown that on $\RP$, the operator $S$ 
coincides with the proximal mapping of the function 
\begin{equation*}
g(t) := \begin{cases}
\tfrac14\ln\big(1+e^{-2f(t)}\big), & \text{if $t\geq 0$;}\\
+\infty, &\text{if $t<0$,}
\end{cases}
\end{equation*}
which is strictly decreasing, strictly convex, and differentiable on 
$\RPP$ with 
$g'=-Q$ while $\inf g(\RR)=0$ is unattained.
}
$S$ that is later responsible for the parameter sequence $(s_n)_{n\geq 0}$: 
\begin{empheq}[box=\mybluebox]{equation}
  \label{e:T}
  S\colon\RP\to\RPP\colon s\mapsto s_+,
  \qquad s_+-s=Q(s_+).
\end{empheq}

\subsection{Convexity of Softplus Compositions}

\label{s:softplus}

We next establish a criterion for the convexity of compositions
with the \emph{softplus function} 
\begin{empheq}[box=\mybluebox]{equation}
  \label{e:dPhi}
  \Phi\colon\RR\to\RR\colon t\mapsto \ln(1+{\color{black}e^{t}}).
\end{empheq}
A direct computation gives 
\begin{equation}
  \label{e:Phi}
  \Phi(t) > 0, \quad
  \Phi'(t) = \frac{{\color{black}e^{t}}}{1+{\color{black}e^{t}}} \in \left]0,1\right[,
  \quad
  \Phi''(t) = \frac{{\color{black}e^{t}}}{(1+{\color{black}e^{t}})^2} > 0.
\end{equation}

\begin{lemma}
  \label{l:Phi}
  Let $I$ be a nondegenerate\footnote{An interval is nondegenerate if 
  its interior is nonempty.} interval in $\RR$, and 
  let $h$ be twice continuously differentiable\footnote{
    A function $h$ is twice continuously differentiable on $I$ 
    if admits a twice continuously differentiable extension to an open interval containing $I$.} 
    on $I$.
  Then
  \begin{equation}
    \Phi\circ h\text{ is convex}
    \Leftrightarrow
    h''+\frac{(h')^2}{1+{\color{black}e^{h}}} \geq 0.
  \end{equation}
\end{lemma}
\begin{proof}
  On the one hand,
  \cref{e:Phi} implies $\Phi''(t)/\Phi'(t) = 1/(1+{\color{black}e^{t}})$.
  On the other hand,
  $(\Phi\circ h)' = \Phi'(h)h'$ and so
  \begin{align*}
    (\Phi\circ h)'' & = (\Phi'(h))'\cdot h'+ \Phi'(h)h''
    = \Phi''(h)h'h' +\Phi'(h)h''                                                  
                    = \Phi'(h)\Big(h'' + \frac{\Phi''(h)}{\Phi'(h)}(h')^2\Big).
  \end{align*}
  Altogether,
  \begin{equation*}
    (\Phi\circ h)'' = \Phi'(h)\Big(h'' + \frac{(h')^2}{1+{\color{black}e^{h}}}\Big),
  \end{equation*}
  and this yields the conclusion.
\end{proof}

\begin{corollary}
  \label{c:Phi}
  Let $I$ be a nondegenerate interval in $\RR$, and 
  let $h$ be twice continuously differentiable on $I$.
  If
  \begin{equation}
    h\leq 0 \;\text{ and }\; h'' + (h')^2/2 \geq 0,
  \end{equation}
  then
  \begin{equation}
    \Phi\circ h = \ln\big(1+{\color{black}e^{h}}\big)\text{ is convex.}
  \end{equation}
\end{corollary}
\begin{proof}
  We have $h\leq 0$
  $\Leftrightarrow$
  ${\color{black}e^{h}}\leq 1$
  $\Leftrightarrow$
  $1+{\color{black}e^{h}}\leq 2$
  $\Leftrightarrow$ $1/(1+{\color{black}e^{h}}) \geq 1/2$.
  It follows that
  \begin{equation}
    h''+\frac{(h')^2}{1+{\color{black}e^{h}}} \geq h'' + \frac{(h')^2}{2} \geq 0,
  \end{equation}
  and the conclusion now follows from \cref{l:Phi}.
\end{proof}

\subsection{\texorpdfstring{A Convexity Property of $f$}{A Convexity Property of f}}

We now return to the function $f$ that satisfies \cref{e:f} and 
show that it induces a convex function via the softplus function $\Phi$ 
from \cref{e:dPhi}: 

\begin{proposition}
  Let $\tau\in\RP$ and set $I := \left[-\tau,+\infty\right[$.
  Then
  \begin{equation}
    \psi_\tau\colon
    I\to\RR\colon x\mapsto
    \ln\Big(1+{\color{black}e^{-f(\tau)-f(\tau+x)+x^2}}\Big)
  \end{equation}
  is a convex function and
  \begin{equation}
    \label{e:Deniz}
    (\forall x\in I)\quad
    \big(1+{\color{black}e^{-2f(\tau)}}\big) {\color{black}e^{-2Q(\tau)x}}
    \leq 1+{\color{black}e^{-f(\tau)-f(\tau+x)+x^2}}.
  \end{equation}
\end{proposition}
\begin{proof}
  Set $h\colon I\to\RR\colon x\mapsto -f(\tau)-f(\tau+x)+x^2$.
By \cref{e:f0}, the function $f$ is twice continuously differentiable
  on an open interval containing $\RP$. Translating that interval by
  $-\tau$ shows that $h$ extends to a twice continuously differentiable
  function on an open interval containing $I$.
  Let $x\in I$, i.e., $x\geq -\tau$.
  Using \cref{ff}\cref{ff:i}, we estimate
  \begin{equation*}
    f(\tau)+f(\tau+x) > \tau^2 + (\tau+x)^2
    = x^2+2\tau(\tau+x) \geq x^2;
  \end{equation*}
  thus,
  \begin{equation*}
    h(x)<0.
  \end{equation*}
  Next, $h'(x) = -f'(\tau+x)+2x = {\color{black}-(f'(\tau+x)-2x)}$
  and $h''(x) = 2-f''(\tau+x)$.
  Since $\tau\geq0$,
  we estimate with the help of \cref{ff}\cref{ff:i} that 
  \begin{equation}
  \label{e:261103e}
    f'(\tau+x)-2x
    =f'(\tau+x)-2(\tau+x)+2\tau
    \geq f'(\tau+x)-2(\tau+x)\geq2.
  \end{equation}
In turn,  \cref{e:f5} and \cref{e:261103e} yield
  \begin{align*}
    h''(x)+\frac{(h'(x))^2}{2}
    &=2-f''(\tau+x)+\frac{(f'(\tau+x)-2x)^2}{2}\\
    &\geq-\frac{(f'(\tau+x)-2(\tau+x))^2}{2}
          +\frac{(f'(\tau+x)-2x)^2}{2}\\
    &\geq0.
  \end{align*}
  Now let $\Phi$ be as in \cref{e:dPhi}. Then $\psi_\tau = \Phi\circ h$.
  The convexity of $\psi_\tau = \Phi\circ h$ now follows from \cref{c:Phi}.
  Next, $\psi_\tau' = (\Phi\circ h)' = \Phi'(h)h'$,
  and so \cref{e:Phi} yields 
  \begin{equation*}
    \psi_\tau'(0) = \Phi'(h(0))h'(0) = \frac{{\color{black}e^{-2f(\tau)}}}{1+{\color{black}e^{-2f(\tau)}}} \cdot (-f'(\tau))
    = -2Q(\tau).
  \end{equation*}
  Hence the subgradient inequality $\psi_\tau(0)+\psi_\tau'(0)x \leq \psi_\tau(x)$
  turns into
  \begin{equation*}
    \ln\big(1+{\color{black}e^{-2f(\tau)}}\big) - 2Q(\tau)x
    \leq
    \ln\big(1+{\color{black}e^{-f(\tau)-f(\tau+x)+x^2}}\big).
  \end{equation*}
  Taking the exponential function on both sides, we
  obtain \cref{e:Deniz}.
\end{proof}

\begin{corollary}
  \label{c:Deniz}
  Let $t\in\RP$, and let $s,s_+$ be as in \cref{p:Q}.
  Then
  \begin{equation}
    \label{e:DDay}
    {\color{black}e^{-(s-t)^2}}
    \leq
    \frac{{\color{black}e^{-(
            s_+-s)^2}}}{1+{\color{black}e^{-2f(s_+)}}}
    \Big({\color{black}e^{-(
          s_+-t)^2}} + {\color{black}e^{-f(s_+)-f(t)}} \Big).
  \end{equation}
\end{corollary}
\begin{proof}
  Inequality \cref{e:Deniz}, with $\tau=s_+$ and $x=t-s_+$,
  reads
  \begin{equation}
    \big(1+{\color{black}e^{-2f(s_+)}}\big) {\color{black}e^{-2Q(s_+)(t-s_+)}}
    \leq 1+{\color{black}e^{-f(s_+)-f(t)+ (t-s_+)^2}},
  \end{equation}
  which, after recalling \cref{s+-s}, turns into
  \begin{equation}
    \big(1+{\color{black}e^{-2f(s_+)}}\big) {\color{black}e^{-2(s_+-s)(t-s_+)}}
    \leq 1+{\color{black}e^{-f(s_+)-f(t)+ (t-s_+)^2}}.
  \end{equation}
  Now divide this inequality
  by $\big(1+{\color{black}e^{-2f(s_+)}}\big) {\color{black}e^{(s_+-s)^2+(s_+-t)^2}}$
  and simplify to obtain exactly \cref{e:DDay}.
\end{proof}

\section{Alternating proximal mappings and their Cesàro means}

\label{s:Phi} 

In this section, we assume that 
\begin{subequations}
\label{eqs:Phi}
\begin{empheq}[box=\mybluebox]{equation}
  \HH \text{ is a real Hilbert space, with inner product $\scal{\cdot}{\cdot}$ and induced norm $\|\cdot\|$.}
\end{empheq}
We also suppose that 
\begin{empheq}[box=\mybluebox]{equation}
  \text{$F,G$ are proper lower semicontinuous convex functions from $\HH$ to
  $\left]-\infty,+\infty\right]$}.
\end{empheq}
We set 
\begin{empheq}[box=\mybluebox]{equation}
  \Phi \colon \HH\times\HH\to\left]-\infty,+\infty\right]\colon 
  (x,y)\mapsto F(x)+G(y)+\tfrac{1}{2}\norm{x-y}^2 
\end{empheq}
and we assume that $\Phi$ has a minimizer, i.e., 
$\inf\Phi(\HH\times\HH)$ is attained: 
\begin{empheq}[box=\mybluebox]{equation}
\label{e:muattained}
  \mu := \min \Phi(\HH\times\HH) \in\RR.
\end{empheq}
Given $y_0\in\HH$, we consider the alternating proximal mapping sequence 
defined by
\begin{empheq}[box=\mybluebox]{equation}
  \label{eq:orbit}
(\forall n\geq 0)\quad
 x_{n+1}=\prox_{F}(y_n),\qquad
 y_{n+1}=\prox_{G}(x_{n+1})
\end{empheq}
as well as the corresponding Cesàro means
\begin{empheq}[box=\mybluebox]{equation}
(\forall n\geq 1)\quad
  \overline{x}_n := \frac{1}{n}\sum_{k=1}^{n} x_k
  \quad\text{and}\quad
  \overline{y}_n := \frac{1}{n}\sum_{k=0}^{n-1} y_k. 
\end{empheq}
\end{subequations}

\begin{theorem}
\label{t:Phi}
Recall the standing assumption \cref{eqs:Phi}. 
Then the sequences defined in \cref{eq:orbit} 
satisfy 
\begin{equation}
\label{e:Phi1}
  \sum_{n=1}^{\infty}n\big(\norm{x_{n+1}-x_n}^2
  + \norm{y_{n+1}-y_n}^2\big)<+\infty, 
\end{equation}
\begin{equation}
\label{e:Phi2}
n(x_{n+1}-x_n)\to 0
\quad\text{and}\quad
n(y_{n+1}-y_n)\to 0, 
\end{equation}
and 
\begin{equation}
\label{e:Phi3}
x_n-\overline{x}_n\to 0
\quad\text{and}\quad 
y_n-\overline{y}_n\to 0. 
\end{equation}
\end{theorem}
\begin{proof}
Set 
\begin{equation}
T := \prox_{G}\circ\prox_{F}.
\end{equation}
Let $(x,y)\in\HH\times\HH$. 
Then $(x,y)$ is a minimizer of $\Phi$ 
$\Leftrightarrow$ 
$0\in\partial F(x)+x-y$ $\land$ 
$0\in\partial G(y)+y-x$
$\Leftrightarrow$
$x=\prox_{F}(y)$ $\land$ 
$y=\prox_{G}(x)$
$\Leftrightarrow$
$y\in\Fix T$ $\land$ $x=\prox_{F}(y)$.
We deduce that $\Fix T$ is nonempty exactly when 
$\argmin\Phi$ is nonempty, which is our assumption \cref{e:muattained}.

Now set 
\begin{equation}
(\forall n\geq 1)\quad 
 \alpha_n:=\Phi(x_{n+1},y_n)-\mu
 \;\;\text{and}\;\;\beta_n:=\Phi(x_{n+1},y_{n+1})-\mu.
\end{equation}
It follows from \cite[Theorem~4.6]{BCR} that 
$(\forall n\geq 1)$ $\alpha_{n}\geq \beta_n\geq \alpha_{n+1}$, 
that $\alpha_n\to 0$ and $\beta_n\to 0$, and that
\begin{equation}
\label{e:261108b}
\sum_{n\geq 1} \alpha_n < \infty \;\;\text{and} \;\;
\sum_{n\geq 1} \beta_n < \infty.
\end{equation}
Fix $k\geq 1$. 
Because $y_{k+1}=\prox_{G}(x_{k+1})$ $\Leftrightarrow$ 
$x_{k+1}-y_{k+1}\in\partial G(y_{k+1})$, we obtain 
$G(y_k) \geq G(y_{k+1}) + \scal{y_k-y_{k+1}}{x_{k+1}-y_{k+1}}$
$\Leftrightarrow$ 
$G(y_k) + \thalb \|y_k-x_{k+1}\|^2 - G(y_{k+1}) - \thalb \|y_{k+1}-x_{k+1}\|^2 
\geq \thalb \|y_k-y_{k+1}\|^2$. 
Hence 
\begin{subequations}
\label{e:261108a}
\begin{align}
\thalb\|y_{k+1}-y_k\|^2
&\leq G(y_k) + \thalb \|y_k-x_{k+1}\|^2 - G(y_{k+1}) - \thalb \|y_{k+1}-x_{k+1}\|^2\\
&= \Phi(x_{k+1},y_k) - \Phi(x_{k+1},y_{k+1})\\
&= \alpha_k - \beta_{k}\\
&\leq \alpha_k - \alpha_{k+1}.
\end{align}
\end{subequations}
Summing $k\cdot\cref{e:261108a}$ over $k\in\{1,\ldots,n\}$ 
and recalling \cref{e:261108b} yields
\begin{align*}
\thalb \sum_{k=1}^{n} k\|y_{k+1}-y_k\|^2
&\leq \sum_{k=1}^{n} k(\alpha_k - \alpha_{k+1})
= \sum_{k=1}^{n} \big(k\alpha_k - (k+1)\alpha_{k+1}+\alpha_{k+1}\big)\\
&= \alpha_1-(n+1)\alpha_{n+1} + \sum_{k=1}^{n} \alpha_{k+1}
=-n\alpha_{n+1} + \sum_{k=1}^{n} \alpha_k\\
&\leq \sum_{k=1}^{\infty} \alpha_k < \infty. 
\end{align*}
It follows that 
\begin{equation}
\label{e:261108c}
\sum_{n=1}^{\infty} n\|y_{n+1}-y_n\|^2 < \infty 
\end{equation}
and also that 
$\sum_{n=1}^{\infty} (n+1)\|y_{n+1}-y_n\|^2 < \infty$. 
Because $\prox_{F}$ is nonexpansive, we deduce that
$\sum_{n=1}^{\infty} (n+1)\|x_{n+2}-x_{n+1}\|^2 < \infty$ which implies 
\begin{equation}
\label{e:261108d}
\sum_{n=1}^{\infty} n\|x_{n+1}-x_n\|^2 < \infty. 
\end{equation}
Combining \cref{e:261108c} and \cref{e:261108d} yields \cref{e:Phi1}.

The mapping $T$ is nonexpansive, because it is the composition of two (firmly) 
nonexpansive mappings.
Hence the sequence $(\|y_{n+1}-y_n\|)_{n\geq 1}$ is decreasing. 
Now let $n\geq 2$. 
The set $\{\lceil n/2\rceil,\ldots,n\}$ has $n-\lceil n/2\rceil+1 = \lfloor n/2\rfloor+1 
\geq n/2$ elements. 
If $k\in\{\lceil n/2\rceil,\ldots,n\}$, then $k\geq n/2$ and 
$\|y_{k+1}-y_k\|\geq \|y_{n+1}-y_n\|$. 
Combining this with \cref{e:261108c}, we obtain 
\begin{align*}
\frac{n^2}{4}\|y_{n+1}-y_n\|^2
&\leq \big(n-\lceil n/2\rceil+1\big)\frac n2\,\|y_{n+1}-y_n\|^2
\leq \sum_{k=\lceil n/2\rceil}^n k \|y_{k+1}-y_k\|^2\\
&\leq \sum_{k=\lceil n/2\rceil}^{\infty}k \|y_{k+1}-y_k\|^2
\to 0 \quad \text{as } n\to\infty.
\end{align*}
In turn, this implies $n(y_{n+1}-y_n)\to 0$. 
Because $\prox_{F}$ is nonexpansive, we deduce that 
$n(x_{n+2}-x_{n+1})\to 0$, which implies $n(x_{n+1}-x_n)\to 0$.
This verifies \cref{e:Phi2}. 

We now turn to the Cesàro means. 
For $n\geq 1$, we have 
\begin{align*}
 y_n-\overline{y}_n
 &=\frac1n\sum_{k=0}^{n-1}(y_n-y_k)
   =\frac1n\sum_{k=0}^{n-1}\sum_{j=k}^{n-1}(y_{j+1}-y_j)\\
 &=\frac1n\sum_{j=0}^{n-1}\sum_{k=0}^{j}(y_{j+1}-y_j)
   =\frac1n\sum_{j=0}^{n-1}(j+1)(y_{j+1}-y_j).
\end{align*}
Therefore,
\begin{equation*}
  \|y_n-\overline{y}_n\|
  \leq\frac1n\sum_{k=0}^{n-1}(k+1)\|y_{k+1}-y_k\|\to0 
  \quad\text{as } n\to\infty,
\end{equation*}
where the last limit follows from the Ces\`aro convergence theorem 
applied to the real sequence $((n+1)\|y_{n+1}-y_n\|)_{n\geq 0}$ 
which tends to zero by \cref{e:Phi2}. 
We've verified that $y_n-\overline{y}_n\to0$. 
The proof of $x_n-\overline{x}_n\to0$ is analogous.
\end{proof}

\begin{remark}
\label{r:DavisYin}
For further properties related to \cref{t:Phi}, we refer the reader to 
\cite{BCR}. 
When $F$ or $G$ is identically equal to $0$, i.e., 
we just iterate one proximal mapping, then 
\cref{e:Phi1} may also be deduced from \cite[Theorem~4.1.3]{ABN}. 
The convergence result \cref{e:Phi2} may also be deduced from 
Davis and Yin's forward--backward residual estimate \cite[Theorem~12 on page~150]{DavisYin}.
Finally, we are unaware of an earlier explicit statement of \cref{e:Phi3} for 
alternating proximal mappings. It may fail for general resolvents rather than 
the proximal mappings considered here: indeed, \cite{BT26} contains an example 
of a resolvent iteration that converges weakly but not strongly to a fixed point but 
where the Cesàro means do converge strongly.
\end{remark}

\section{\texorpdfstring{Properties of $\bH$ and $\bK$}{Properties of H and K}}
\label{s:sets}

We recall from \cref{e:bH} that $\bH = \{0\}\times X$. 
Let us record the description of the polar cone of $\bK$,
which follows immediately from the definition of $\bK$ in \cref{e:bK}:
\begin{equation}
  \label{e:bK-}
  \bK^\ominus
  = \menge{(\beta,y)\in \mathbb R\times X}%
  {(\forall t\in\RP)\; \beta{\color{black}e^{-f(t)}} + \scal{y}{u(t)} \leq 0}.
\end{equation}

\begin{proposition}
\label{p:HK0}
  We have
  \begin{equation}
    \bH\cap \bK = \{\bzero\}.
  \end{equation}
\end{proposition}
\begin{proof}
  It is clear that $\bzero = (0,0)\in \bH\cap \bK$.
  Conversely, let $(0,x)\in \bH\cap \bK$, and set
  \begin{equation*}
    U := \cspan u(\RP) = \cspan\{u(t)\}_{t\in\RP}.
  \end{equation*}
  It's clear that 
  \begin{equation}
  \label{e:261105a}
  x\in U
  \end{equation}
  and we must show that $x=0$.
  To this end, fix $r\in\RP$.
  Recalling \cref{e:Kurve}, we obtain
  $(\forall t\in\RP)$
  $\tscal{(0,-u(r))}{({\color{black}e^{-f(t)}},u(t))} =
    \scal{-u(r)}{u(t)} = -{\color{black}e^{-(t-r)^2}} \leq 0$;
  hence, by \cref{e:bK-}, we have $(0,-u(r))\in \bK^\ominus$.
  Because $(0,x)\in \bK$, we obtain
  $\scal{(0,x)}{(0,-u(r))} = -\scal{x}{u(r)} \leq 0$;
  thus,
  \begin{equation}
    \label{e:Vati1}
    \scal{x}{u(r)} \geq 0.
  \end{equation}
  Now pick any $R>r$ and $\rho>0$, and set
  \begin{equation*}
    v := u(r) - \rho u(R).
  \end{equation*}
  By \cref{e:Kurve}, we have
  \begin{align*}
    \scal{v}{u(t)} 
    &= \scal{u(r)-\rho u(R)}{u(t)}
    = {\color{black}e^{-(t-r)^2}} - \rho {\color{black}e^{-(t-R)^2}}
    = e^{-(t-R)^2}\Big( \frac{e^{-(t-r)^2}}{e^{-(t-R)^2}} - \rho \Big).
  \end{align*}
  Because $R>r$, we deduce that
  \begin{equation*}
    \frac{{\color{black}e^{-(t-r)^2}}}{{\color{black}e^{-(t-R)^2}}} = {\color{black}e^{-(t-r)^2+(t-R)^2}}
    = {\color{black}e^{R^2-r^2-2t(R-r)}}\to 0 \quad\text{as }t\to\infty.
  \end{equation*}
  Because $\rho>0$, it now follows that
  $\scal{v}{u(t)} < 0$ for all sufficiently large $t$.
  On the remaining compact interval the function
  $t\mapsto{\color{black}e^{f(t)}}\scal{v}{u(t)}$ is continuous and hence bounded above.
  Thus
  \begin{equation*}
    \sigma := \sup_{t\in\RP} {\color{black}e^{f(t)}}\scal{v}{u(t)} < \infty.
  \end{equation*}
  We now obtain
  $(\forall t\in\RP)$
  $\scal{v}{u(t)} \leq \sigma {\color{black}e^{-f(t)}}$;
  thus, by \cref{e:bK-}, we have
  \begin{equation*}
    (-\sigma, v) \in \bK^\ominus.
  \end{equation*}
  Recalling that $(0,x)\in \bK$, we get
  $0 \geq \scal{(0,x)}{(-\sigma,v)}
    = \scal{x}{v}
    = \scal{x}{u(r)} - \rho \scal{x}{u(R)}$.
  Letting now $\rho\to 0^+$, we deduce that
  \begin{equation}
    \label{e:Vati2}
    \scal{x}{u(r)} \leq 0.
  \end{equation}
  Combining \cref{e:Vati1} with \cref{e:Vati2}, we obtain $\scal{x}{u(r)}=0$.
  Finally, since $r\in\RP$ was arbitrary, we obtain
  $x\in U^\perp$. 
  In view of \cref{e:261105a}, we conclude that $x=0$.
\end{proof}

We now turn to the projections onto $\bH$ and $\bK$.
Clearly,
\begin{equation}
  \label{PH}
  \Pj{\bH}\colon \bX\to \bX\colon (\alpha,x)\mapsto (0,x).
\end{equation}
While we don't have a complete description of $\Pj{\bK}$,
the following will be sufficient for our purposes:

\begin{proposition}
  Let $\rho \in\RP$ and $s\in\RP$.
  Suppose that $s_+$ is as in \cref{p:Q}.
  Then
  \begin{subequations}
    \label{PK}
    \begin{equation}
      \Pj{\bK}\rho\big(0,u(s)\big) = \rho_+\big({\color{black}e^{-f(s_+)}}, u(s_+)\big) = \rho_+\bw(s_+),
    \end{equation}
    where
    \begin{equation}
      \label{rho+}
      \rho_+ := \rho\frac{{\color{black}e^{-(s_+-s)^2}}}{1+{\color{black}e^{-2f(s_+)}}}.
    \end{equation}
  \end{subequations}
\end{proposition}
\begin{proof}
  The result is clear if $\rho=0$, so we assume that $\rho>0$.
  Note that
  \begin{equation}
    \label{e:klu1}
    \rho_+\bw(s_+) \in \bK.
  \end{equation}
  Let $t\in\RP$.
  Using \cref{e:Kurve}, 
   we obtain
  \begin{subequations}
    \label{e:gut}
    \begin{align}
       & \hspace{-0.5cm}\scal{\rho(0,u(s)) - \rho_+\bw(s_+)}{\bw(t)}                                                              \\
       & = \scal{(0,\rho u(s)) - \big(\rho_+{\color{black}e^{-f(s_+)}}, \rho_+u(s_+)\big)}{\big({\color{black}e^{-f(t)}},u(t)\big)} \\
       & = \scal{\big(-\rho_+{\color{black}e^{-f(s_+)}},\rho u(s) - \rho_+u(s_+)\big)}{\big({\color{black}e^{-f(t)}},u(t)\big)}     \\
       & = -\rho_+{\color{black}e^{-f(s_+)-f(t)}} + \rho\scal{u(s)}{u(t)} - \rho_+\scal{u(s_+)}{u(t)}                              \\
       & = -\rho_+{\color{black}e^{-f(s_+)-f(t)}} + \rho{\color{black}e^{-(s-t)^2}} - \rho_+{\color{black}e^{-(s_+-t)^2}}.
    \end{align}
  \end{subequations}
  Substituting $t=s_+$ into \cref{e:gut} and recalling \cref{rho+}, we get
  \begin{align*}
    \scal{\rho(0,u(s)) - \rho_+\bw(s_+)}{\bw(s_+)}
     & = -\rho_+{\color{black}e^{-2f(s_+  )}} + \rho{\color{black}e^{-(s_+-s)^2}} - \rho_+ 
      = 0;
  \end{align*}
  consequently,
  \begin{equation}
    \label{e:klu2}
    \scal{\rho(0,u(s)) - \rho_+\bw(s_+)}{\rho_+\bw(s_+)} = 0.
  \end{equation}
  Continuing from \cref{e:gut},
  recalling \cref{rho+}, and using
  \cref{e:DDay} for the following inequality, we estimate
  \begin{align*}
     \scal{\rho(0,u(s)) - \rho_+\bw(s_+)}{\bw(t)}                                                                                                                         
     & = -\rho_+{\color{black}e^{-f(s_+)-f(t)}} + \rho{\color{black}e^{-(s-t)^2}} - \rho_+{\color{black}e^{-(s_+-t)^2}}                                                                       \\
     & =\rho\Big({\color{black}e^{-(s-t)^2}} - \frac{\rho_+}{\rho}\big({\color{black}e^{-f(s_+)-f(t)}} +{\color{black}e^{-(s_+-t)^2}}\big)\Big)                                               \\
     & =\rho\Big({\color{black}e^{-(s-t)^2}} - \frac{{\color{black}e^{-(s_+-s)^2}}}{1+{\color{black}e^{-2f(s_+)}}}\big({\color{black}e^{-f(s_+)-f(t)}} +{\color{black}e^{-(s_+-t)^2}}\big)\Big) \\
     & \leq 0.
  \end{align*}
  Because $t\in\RP$ was arbitrary, we learn from
  the last inequality that
  \begin{equation}
    \label{e:klu3}
    \rho\big(0,u(s)\big) - \rho_+\bw(s_+)\in \bK^\ominus.
  \end{equation}
  Combining \cref{e:klu1}, \cref{e:klu2}, \cref{e:klu3},
  with \cite[Proposition~6.28]{BC2017}, we deduce that $\rho_+\bw(s_+) = \Pj{\bK}\rho(0,u(s))$, and we're done.
\end{proof}

\begin{remark}
  In view of \cref{e:T}, we can also write \cref{PK} as
  \begin{equation}
    \label{PKT}
    \Pj{\bK}\rho\big(0,u(s)\big) = \rho_+\bw\big(S(s)\big),
    \;\;\text{where }\;\;
    \rho_+ := \rho\frac{{\color{black}e^{-(S(s)-s)^2}}}{1+{\color{black}e^{-2f(S(s))}}}.
  \end{equation}
\end{remark}

\section{The Sequence of Alternating Projections}
\label{s:seq}

We now set up the sequence of alternating projections.
The starting point is (recall \cref{e:bw})
\begin{empheq}[box=\mybluebox]{equation}
  \label{e:y0}
  \by_0 := \rho_0\bw(s_0) =
  \rho_0\big({\color{black}e^{-f(s_0)}}, u(s_0)\big)\in\bK,
  \quad\text{where $\rho_0>0$ and $s_0 \geq 0$.}
\end{empheq}
Given $\by_n = \rho_n\bw(s_n)\in\bK$, we update
by projecting first onto $\bH$ and then onto $\bK$:
\begin{empheq}[box=\mybluebox]{equation}
  \label{e:xnyn}
  \bx_{n+1} := \Pj{\bH}(\by_n)
  \quad\text{and}\quad
  \by_{n+1} := \Pj{\bK}(\bx_{n+1}) = \Pj{\bK}\Pj{\bH}(\by_n).
\end{empheq}
The two sequences $(\bx_n)_{n\geq 1}$ and $(\by_n)_{n\geq 0}$
together form the sequence of alternating projections
$(\by_0,\bx_1,\by_1,\bx_2,\by_2,\ldots)$.
These sequences, which live in $\bX$,
are fully described by two scalar sequences
$(\rho_n)_{n\geq 0}$ and $(s_n)_{n\geq 0}$, generated
iteratively by updating given $s_n$ and $\rho_n$ via
\begin{empheq}[box=\mybluebox]{equation}
  \label{srho}
  s_{n+1} := S(s_n) 
  \quad\text{and}\quad
  \rho_{n+1} := \rho_n\frac{{\color{black}e^{-(s_{n+1}-s_n)^2}}}{1+{\color{black}e^{-2f(s_{n+1})}}},
\end{empheq}
where $S$ is defined in \cref{e:T}.

\begin{proposition}
  The alternating projections sequences
  are exactly described by the formulas
  \begin{equation}
    \label{e:xnyn'}
    (\forall n\geq 1)\quad
    \bx_n = \rho_{n-1}\big(0,u(s_{n-1})\big)
    \;\;\text{and}\;\;
    \by_n = \rho_n\bw(s_n).
  \end{equation}
\end{proposition}
\begin{proof}
  This follows by induction from
  \cref{e:xnyn}, \cref{PH}, \cref{PKT}, and \cref{srho}.
\end{proof}

Thus to understand the behaviour of the alternating projections
sequences, we must understand the behaviour of the scalar sequences $(s_n)_{n\geq 0}$ and $(\rho_n)_{n\geq 0}$.

\begin{proposition}[properties of $(s_n)_{n\geq 0}$]
  \label{p:sn}
  The sequence $(s_n)_{n\geq 0}$ satisfies the following:
  \begin{subequations}
  \begin{equation}
    \label{e:sn1}
    0 \leq s_0 < s_1 < \cdots < s_n < s_{n+1}\to \infty,
  \end{equation}
  \begin{equation}
    \label{e:sn2}
    \sum_{n=0}^\infty {e^{-2f(s_{n})}} <\infty,
  \end{equation}
  \begin{equation}
    \label{e:sn3}
    (\forall m\in\{0,1,\ldots\})\quad
    \sum_{k=m}^\infty (s_{k+1}-s_k)^2 < \tfrac{1}{4} {e^{-2f(s_m)}},
  \end{equation}
  \begin{equation}
    \label{e:sn4}
    \sum_{n=0}^\infty (s_{n+1}-s_n)^2 < \infty.
  \end{equation}
\end{subequations}
\end{proposition}
\begin{proof}
  Strict increase of $(s_n)_{n\geq 0}$ is clear from \cref{srho}, \cref{e:T},
  and \cref{svss+}.
  Set $s_\infty := \lim_{n\to\infty} s_n \in \left]0,+\infty\right]$.
  Suppose to the contrary that $s_\infty < +\infty$.
  By \cref{s+-s} and \cref{p:Q}\cref{p:Q1}\&\cref{p:Q2}, we have
  \begin{equation*}
    0\leftarrow s_{n+1}-s_n = Q(s_{n+1}) \to Q(s_\infty)>0,
  \end{equation*}
  which is absurd. We've verified \cref{e:sn1}.

  Next, using 
  \cref{Es+}, 
  \cref{e:sn1},
  \cref{ff}\ref{ff:i}, 
  and \cref{e:f2},
  we obtain
  \begin{align*}
    \sum_{n=0}^\infty e^{-2f(s_{n+1})}
     & <\frac{5}{2}\sum_{n=0}^\infty \int_{s_n}^{s_{n+1}} \frac{dt}{f'(t)}
    = \frac{5}{2}\int_{s_0}^{\infty} \frac{dt}{f'(t)}
    \leq
    \frac{5}{2}\int_{0}^{\infty} \frac{dt}{f'(t)}<\infty.
  \end{align*}
  This and the fact that
  ${\color{black}e^{-2f(s_0)}}<\infty$ yield \cref{e:sn2}.

By \cref{Es-Es+}, for every $k\geq m\geq0$, we get
\begin{equation*}
  4(s_{k+1}-s_k)^2
  <e^{-2f(s_k)}-e^{-2f(s_{k+1})}.
\end{equation*}
This sums to 
\cref{e:sn3}, which in turn implies \cref{e:sn4}.
\end{proof}

\begin{proposition}[properties of $(\rho_n)_{n\geq0}$]
\label{p:rhon}
The sequence $(\rho_n)_{n\geq0}$ is strictly decreasing and converges
 to a positive limit $\rho_\infty$:
\begin{equation}
\rho_0 > \rho_1 > \cdots > \rho_n > \rho_{n+1} \to \rho_\infty > 0.
\end{equation}
Moreover, 
\begin{subequations}
\begin{align}
\label{e:rholower}
  \rho_\infty
  &>\rho_0\exp\left(
    -\frac14e^{-2f(s_0)}
    -\frac52\int_{s_0}^{\infty}\frac{dt}{f'(t)}
    \right)\\
\label{e:rholoweruniform}
  &\geq\rho_0\exp\left(-\frac1{16}-\frac52
  \int_0^\infty\frac{dt}{f'(t)}\right)\\
  &>0.
\end{align}
\end{subequations}
\end{proposition}
\begin{proof}
Because $\rho_0>0$, \cref{srho} shows that $(\rho_n)_{n\geq0}$ is
strictly decreasing and positive. Hence it converges to some
$\rho_\infty\geq0$.
By \cref{e:sn3} with $m=0$, we have
\begin{equation*}
  \sum_{k=0}^{\infty}(s_{k+1}-s_k)^2<\frac14e^{-2f(s_0)}. 
\end{equation*}
Moreover, summing \cref{Es+} over the successive intervals
$[s_k,s_{k+1}]$, and using \cref{e:sn1} and \cref{e:f2}, gives
\begin{equation*}
  \sum_{k=0}^{\infty}e^{-2f(s_{k+1})}
  <\frac52\sum_{k=0}^{\infty}
       \int_{s_k}^{s_{k+1}}\frac{dt}{f'(t)}
  =\frac52 \int_{s_0}^{\infty}\frac{dt}{f'(t)}<\infty.
\end{equation*}
Since $1\leq 1+t\leq e^t$ 
$\Leftrightarrow$ $0\leq\ln(1+t)\leq t$ for every $t\geq0$, it follows that
\begin{equation}
\label{e:sigma}
  \sigma:=\sum_{k=0}^{\infty}
       \big((s_{k+1}-s_k)^2+\ln(1+e^{-2f(s_{k+1})})\big)
  <\frac14 e^{-2f(s_0)}+\frac52 \int_{s_0}^\infty \frac{dt}{f'(t)}<\infty.
\end{equation}
For every integer $n\geq1$, the recurrence \cref{srho} gives
\begin{align*}
  \rho_n
  &=\rho_0\prod_{k=0}^{n-1}
       \frac{e^{-(s_{k+1}-s_k)^2}}{1+e^{-2f(s_{k+1})}}
  =\rho_0\exp\left(
       -\sum_{k=0}^{n-1}
          \big((s_{k+1}-s_k)^2+\ln(1+e^{-2f(s_{k+1})})\big)
       \right).
\end{align*}
Letting $n\to\infty$ and recalling \cref{e:sigma}, we thus obtain
\begin{equation*}
  \rho_\infty=\rho_0e^{-\sigma}
  >\rho_0\exp\left(-\frac14e^{-2f(s_0)}-\frac52\int_{s_0}^{\infty}\frac{dt}{f'(t)}\right)>0,
\end{equation*}
which proves \cref{e:rholower}.
Finally, \cref{ff}\ref{ff:iii} gives $e^{-2f(s_0)}\leq1/4$, while
$\int_{s_0}^{\infty}\frac{dt}{f'(t)}\leq \int_{0}^{\infty}\frac{dt}{f'(t)}$ follows from $s_0\geq0$ and $f'>0$ 
(see \cref{ff}\cref{ff:i}).
These two inequalities yield \cref{e:rholoweruniform}.
\end{proof}

\section{Weak Convergence only!} 

\label{s:main}

Recall that we defined the alternating projections sequences $(\bx_n)_{n\geq 1}$ and $(\by_n)_{n\geq 0}$ in \cref{e:xnyn}, 
and we let $\rho_\infty$ denote the positive limit of $(\rho_n)_{n\geq 0}$ (see \cref{p:rhon}).

\begin{theorem}[the alternating projections sequences converge weakly but not strongly]
  \label{t:map}
  We have
  \begin{equation}
    \label{e:map1}
    \|\by_0\|\geq \|\bx_1\|\geq \|\by_1\|\geq \cdots \geq
    \|\bx_n\|\geq \|\by_n\|\geq \|\bx_{n+1}\| \geq \|\by_{n+1}\|
    \geq \cdots,
  \end{equation}
  $\rho_\infty>0$, 
  and
  \begin{equation}
    \label{e:map2}
    \bx_n \weakly \bzero,\;\;
    \by_n\weakly \bzero,\;\;
    \bx_n-\by_n\to \bzero,\;\;
    \by_n-\bx_{n+1}\to\bzero,
    \;\;\text{but}\;\;
    \|\bx_n\|\to\rho_\infty
    \;\;\text{and}\;\;
    \|\by_n\|\to\rho_\infty.
  \end{equation}
\end{theorem}
\begin{proof}
  Because $s_n\to\infty$ by \cref{e:sn1}, we have
  $u(s_n)\weakly 0$ by \cref{e:Kweak}.
  On the other hand, $\rho_n\to \rho_\infty>0$ by \cref{p:rhon}.
Also, \cref{ff}\ref{ff:i} and \cref{e:sn1} give
  $f(s_n)\to\infty$, and hence $e^{-f(s_n)}\to0$.
  Combining these with \cref{e:xnyn'}, we obtain altogether
  \begin{equation*}
    \bx_n = \rho_{n-1}\big(0,u(s_{n-1})\big) \weakly \rho_\infty(0,0) =\bzero
  \end{equation*}
  and (recall \cref{e:bw})
  \begin{equation*}
    \by_n = \rho_{n}\bw(s_n) = \rho_n\big({\color{black}e^{-f(s_n)}},u(s_n)\big)
    \weakly \rho_\infty(0,0) =\bzero
  \end{equation*}
  as well\footnote{This also follows from
    the general theory of alternating projections.}.
  Moreover, \cref{e:xnyn'}, \cref{e:Kurve},
  and \cref{p:rhon} imply
  \begin{equation}
  \label{e:66}
    \|\bx_n\| = \rho_{n-1}\|u(s_{n-1})\| = \rho_{n-1} \to \rho_\infty > 0.
  \end{equation}
  Because $\Pj{\bH}$ and $\Pj{\bK}$ are nonexpansive and
  $\bzero\in\bH\cap\bK = \Fix(\Pj{\bH}) \cap \Fix(\Pj{\bK})$,
  we have
  \begin{align*}
    \|\by_0-\bzero\|
     & \geq  \|\Pj{\bH}(\by_0) - \Pj{\bH}(\bzero)\|
    = \|\bx_1-\bzero\|
    \geq \|\Pj{\bK}(\bx_1) - \Pj{\bK}(\bzero)\|
    = \|\by_1-\bzero\|                              \\
     & \geq \|\bx_2-\bzero\|
    \geq \cdots;
  \end{align*}
therefore, by \cref{e:66}, we have 
  \begin{equation*}
    \|\by_0\|\geq \|\bx_1\|\geq \|\by_1\|
    \geq \cdots \geq \|\bx_n\|\geq \|\by_n\|\geq \|\bx_{n+1}\|
    \to \rho_\infty.
  \end{equation*}
  Because $\Pj{\bH}$ and $\Pj{\bK}$ are \emph{firmly} nonexpansive and
  $\bzero\in\bH\cap\bK = \Fix(\Pj{\bH}) \cap \Fix(\Pj{\bK})$,
  we now get for $n\geq 0$ that
  \begin{align*}
    \|\by_{n}-\bx_{n+1}\|^2
     & = \|(\Id-\Pj{\bH})(\by_n) -
    (\Id-\Pj{\bH})(\bzero)\|^2                                      \\
     & \leq \|\by_n-\bzero\|^2 - \|\Pj{\bH}\by_n-\Pj{\bH}\bzero\|^2 \\
     & = \|\by_n\|^2 - \|\bx_{n+1}\|^2                              \\
     & \to \rho_\infty^2 - \rho_\infty^2 = 0;
  \end{align*}
similarly, for $n\geq 1$, we have
    $\|\bx_{n}-\by_{n}\|^2
     \leq \|\bx_{n}\|^2 - \|\by_{n}\|^2 \to 0$.
\end{proof}

\begin{remark}[comparison with Hundal's construction]
\label{r:Hundal}
Both Hundal's construction and ours lift a curve of unit vectors by
a positive height and take the closed convex cone generated by the
lifted curve. In both cases, a supporting inequality is needed to
verify the projection onto the entire cone. Hundal already gives
scalar projection equations and estimates within a coordinate block;
see \cite[Theorem~14]{Hundal}. He also obtains a positive numerical
lower bound for the norms of his selected orbit; see
\cite[Theorem~22]{Hundal}.

The distinction here is the single key identity
$\scal{u(s)}{u(t)}=e^{-(s-t)^2}$, valid for every pair of parameters.
It avoids the separate coordinate-plane transition arguments in
Hundal's construction. Moreover, the same construction applies
to every $f$ satisfying \cref{e:f}, and \cref{p:rhon} gives lower
bounds for the limiting radius directly in terms of $f$ and the
initial data.

We do not have a formula for the projection onto $\bK$ at every
point of $\bX$; \cref{PK} covers only the structured inputs needed for
the prescribed family of orbits. There is also a difference in the
generating curves: Hundal's generator has at most three nonzero
coordinates at each parameter, whereas the displayed Gaussian
realization generally has infinitely many. Finally, the
related constructions of Matou\v{s}kov\'a--Reich \cite{MR2003} and
Kopeck\'a \cite{Kopecka2008} also include disjoint-set variants,
and Kopeck\'a obtains a cone with nonempty interior. 
\end{remark}

We now turn to the corresponding Ces\`aro means which are defined by
\begin{empheq}[box=\mybluebox]{equation}
  \overline{\bx}_n := \frac{1}{n}\sum_{k=1}^{n} \bx_k
  \quad\text{and}\quad
  \overline{\by}_n := \frac{1}{n}\sum_{k=0}^{n-1} \by_k
\end{empheq}

\begin{corollary}[the Ces\`aro means do not converge strongly]
\label{t:cesaro}
We have
\begin{equation}
\label{e:68}
  \bx_n-\overline{\bx}_n\to0
  \quad\text{and}\quad
  \by_n-\overline{\by}_n\to0.
\end{equation}
Consequently,
\begin{equation*}
  \overline{\bx}_n\weakly\bzero,
  \quad\overline{\by}_n\weakly\bzero,
  \quad\|\overline{\bx}_n\|\to\rho_\infty>0,
  \quad\|\overline{\by}_n\|\to\rho_\infty>0.
\end{equation*}
The lower bounds in \cref{e:rholower} and \cref{e:rholoweruniform} therefore also
bound the limiting norms of both sequences of Ces\`aro means.
\end{corollary}
\begin{proof}
Applying \cref{t:Phi} 
(with $F=\iota_\bH$, $G=\iota_\bK$, and $y_0=\by_0$), we obtain \cref{e:68}. 

Combining \cref{e:68} with \cref{t:map},
we obtain the weak convergence of both sequences of means.
Moreover, the reverse triangle inequality gives
\begin{equation*}
  \big|\|\bx_n\|-\|\overline{\bx}_n\|\big|
  \leq\|\bx_n-\overline{\bx}_n\|\to0,
\end{equation*}
and likewise for $(\overline{\by}_n)_{n\geq1}$.
The norm limits therefore follow from \cref{t:map}, and their lower
bounds are those established in \cref{p:rhon}.
\end{proof}

\begin{remark}[Hundal's original example revisited]
The weak-but-not-strong convergence of the alternating projections sequence in 
Hundal's seminal example was established in \cite[Theorem~22]{Hundal}. 
Combining this with \cref{t:Phi}, we now see that the Ces\`aro means of Hundal's original example also 
converge weakly to the origin, but not strongly.
\end{remark}

\begin{remark}[Ces\`aro means and infinite-dimensional optimization]
\label{r:optimization}
  Let us point out that \cref{t:map} and \cref{t:cesaro} can
  be interpreted in different ways.
  First, the alternating projections operator is also
  a projected gradient operator because of
  \begin{equation}
    \Pj{\bK}\Pj{\bH} = \Pj{\bK}\circ \big(\bId- \thalb\nabla d_{\bH}^2\big);
  \end{equation}
  so $(\by_n)_{n\geq 0}$ is a projected-gradient iteration.

  Now set $\bK_1 := \bK$ and let
  $\bK_2 = (2\Pj{\bH}-\bId)(\bK_1)$ be the reflection of $\bK_1$ across $\bH$.
  Set $\bz_1 := \bx_1$. Similar to
  \cite[Theorem~5.1]{BMR}, we have
  \begin{equation}
    \bz_{n+1} :=
    \big(\thalb\Pj{\bK_1} + \thalb\Pj{\bK_2}\big)(\bz_n) = \bx_{n+1}.
  \end{equation}
  Moreover, the operator
  \begin{equation}
    \thalb\Pj{\bK_1} + \thalb\Pj{\bK_2}
  \end{equation}
  is not only a \emph{proximal mapping} (see \cite[Theorem~6.1 and Corollary~7.1]{BMR}),
  but also the \emph{gradient descent operator} 
  (see \cite[Remark~4.8(ii)]{BCR})
  \begin{equation}
    \bId -\nabla\big(\tfrac{1}{4}d_{\bK_1}^2+\tfrac{1}{4}d_{\bK_2}^2 \big)
  \end{equation}
  for the function $\tfrac{1}{4}d_{\bK_1}^2+\tfrac{1}{4}d_{\bK_2}^2$.
  So $(\bx_n)_{n\geq 1}$ arises as a \emph{proximal point iteration}
  and as a \emph{gradient descent iteration}.

In view of \cref{t:cesaro}, the Ces\`aro means of these particular
optimization iterations converge weakly to $\bzero$, and their norms
converge to the same positive number $\rho_\infty$.
The bounds in \cref{e:rholower} and \cref{e:rholoweruniform} apply to this limit.
These are consequences of the present construction through the
established relations above.
\end{remark}

\section*{Acknowledgments}
\small
The research of HHB was partially supported by a Discovery Grant
of the Natural Sciences and Engineering Research Council of
Canada.

\section*{AI statement}
\small
The authors used GPT-6 Astra to assist with manuscript preparation and
to test the mathematical arguments. It proposed \cref{p:cutesuff},
assisted with the lower bounds in
\cref{p:rhon}, and
brought relevant references to the authors' attention.
The main ideas for the construction are due to the authors who also 
assume full responsibility for the content of the paper.

\end{document}